\documentclass{article}
\usepackage[english]{babel} 
\usepackage{latexsym,amssymb,amsmath}
\usepackage[cp1251]{inputenc}
\usepackage{amsfonts,amssymb}
\usepackage{graphicx,graphics,hhline}
\usepackage{euscript}
\newcounter{theorem} 
\newcounter{lemma} 

\renewcommand{\thetheorem}{\arabic{theorem}}
\renewcommand{\thelemma}{\arabic{lemma}}
\newcommand{\theor}{\par\refstepcounter{theorem}%
 {\bf Theorem \thetheorem .}\,\,}

\newcommand{\lem}{\par\refstepcounter{lemma}%
{\bf Lemma \thelemma .}\,\,}
\,
\,

\begin{document}
\Large

\begin{center}
\textbf{Constructive description of monogenic functions in
a finite-dimensional commutative associative algebra}
\vskip5mm
 \textbf{V.~S.~Shpakivskyi}
\end{center}
\vskip5mm
\large

\textbf{Abstract.} 
Let $\mathbb{A}_n^m$ be an arbitrary $n$-dimensional 
commutative associative algebra 
over the field of complex numbers with $m$ idempotents. Let
$e_1=1,e_2,e_3$ be elements of $\mathbb{A}_n^m$ which are linearly
independent over the field of real numbers. We consider monogenic
(i.~e. continuous and differentiable in the sense of Gateaux)
functions of the variable $xe_1+ye_2+ze_3$\,, where $x,y,z$ are
real, and obtain a constructive description of all mentioned
functions by means of holomorphic functions of complex variables.
It follows from this description that monogenic functions have
Gateaux derivatives of all orders.

\vskip5mm
\Large

\section{Introduction.} An effectiveness of the analytic function methods
in the complex
plane for researching plane potential fields inspires
mathematicians to develop analogous methods for spatial fields.

Apparently, W.~Hamilton (1843) made the first attempts to
construct an algebra associated with the three-dimensional Laplace
equation\medskip
\begin{equation} \label{Lap3}
 \Delta_3 u(x,y,z):=
\left(\frac{{\partial}^2}{{\partial x}^2}+
 \frac{{\partial}^2}{{\partial y}^2}+
\frac{{\partial}^2}{{\partial z}^2}\right)u(x,y,z)
 =0\,
\end{equation}\medskip
meaning that components of hypercomplex functions satisfy the
equation
(\ref{Lap3}). 
He constructed an algebra of noncommutative quaternions over the
field of real numbers $\mathbb{R}$ and made a base for developing
the hypercomplex analysis.

C.~Segre \cite{Segre} constructed an algebra of commutative
quaternions over the field $\mathbb{R}$ that can be considered as
a two-dimensional commutative semi-simple algebra of bicomplex
numbers over the field of complex numbers $\mathbb{C}$. F.~Ringleb
\cite{Ringleb} and J.~Riley \cite{Riley} obtained a constructive
description of analytic function of a bicomplex variable, namely,
they proved that such an analytic function can be constructed with
an use of two holomorphic functions of complex variables.

Furthermore, F.~Ringleb \cite{Ringleb} considered an arbitrary
finite-dimensional associative (commutative or not)  semi-simple
algebra over the field 
$\mathbb{R}$. For analytic functions which maps the mentioned
algebra onto itself, he obtained a constructive description by
means of real and complex analytic functions.

A relation between spatial potential fields and analytic functions
given in commutative algebras was established by P.~W.~Ketchum
\cite{Ketchum-28} who shown that every analytic function
$\Phi(\zeta)$ of the variable $\zeta=xe_1+ye_2+ze_3$ satisfies the
equation (\ref{Lap3}) in the case where the elements $e_1, e_2,
e_3$ of a commutative algebra satisfy the condition
\begin{equation}\label{garmonichnyj_bazys-ogljad}
     e_1^2+e_2^2+e_3^2=0\,,
\end{equation}
because
\begin{equation}\label{garm}
\frac{{\partial}^{2}\Phi}{{\partial x}^{2}}+
\frac{{\partial}^{2}\Phi}{{\partial y}^{2}}+
\frac{{\partial}^{2}\Phi}{{\partial z}^{2}}\equiv{\Phi}''(\zeta) \
(e_1^2+e_2^2+e_3^2)=0\,,\medskip
\end{equation}
where $\Phi'':=(\Phi')'$ and $\Phi'(\zeta)$ is defined by the
equality $d\Phi=\Phi'(\zeta)d\zeta$.

We say that a commutative associative algebra $\mathbb A$ is {\it
harmonic\/} (cf.
\cite{Ketchum-28,Mel'nichenko75,Plaksa}) if in
$\mathbb A$ there exists a triad of linearly independent vectors
$\{e_1,e_2,e_3\}$ satisfying the equality
\eqref{garmonichnyj_bazys-ogljad} with $e_k^2\ne 0$ for
$k=1,2,3$. We say also that such a triad $\{e_1,e_2,e_3\}$ is {\it
harmonic}.

P.~W.~Ketchum \cite{Ketchum-28} considered the C.~Segre algebra
of quaternions \cite{Segre} as an example of harmonic algebra.

I.~P.~Mel'nichenko \cite{Mel'nichenko75}  noticed that doubly
differentiable in the sense of Gateaux functions form the largest
algebra of functions $\Phi$ satisfying identically the equalities
\eqref{garm}, where $\Phi''$ is the Gateaux second derivative of
function $\Phi$. He proved that there exist exactly $3$
three-dimensional harmonic algebras with unit over the field
$\mathbb C$ (see
\cite{Mel'nichenko75,Melnichenko03,Plaksa}).

Constructive descriptions of monogenic (i.~e. continuous and
differentiable in the sense of Gateaux) functions taking values in
the mentioned three-dimensional harmonic algebras by means three
corresponding holomorphic functions of the complex variable are
obtained in the papers \cite{Pl-Shp1,Pl-Pukh,Pukh}. Such
descriptions make it possible to prove the infinite
differentiability in the sense of Gateaux of monogenic functions
and integral theorems for these functions that are analogous to
classical theorems of the complex analysis (see, e.~g.,
\cite{Pl-Shp3,Plaksa12}).

Furthermore, constructive descriptions of monogenic functions
taking values in special $n$-dimensional 
commutative algebras by means $n$ holomorphic functions of complex
variables are obtained in the papers \cite{Pl-Shp-Algeria,
Pl-Pukh-Analele}.


In this paper we obtaine a constructive description of monogenic
functions taking values in an arbitrary finite-dimensional
commutative associative algebra with unit by means of holomorphic
functions of complex variables.
It follows from this description
  that the monogenic functions have the Gateaux derivatives of all orders.

\section{The algebra $\mathbb{A}_n^m$.}

Let $\mathbb{N}$ be the set of natural numbers.
We fix the numbers $m,n\in\mathbb{N}$ such that $m\leq n$.
Let $\mathbb{A}_n^m$ be an arbitrary commutative associative algebra with
 unit over the field of complex number $\mathbb{C}$.
   E.~Cartan \cite[p.~33]{Cartan}
  proved that there exist a basis $\{I_k\}_{k=1}^{n}$ in $\mathbb{A}_n^m$
  satisfying the following multiplication rules:
  \vskip3mm
1.  \,\,  $\forall\,r,s\in[1,m]\cap\mathbb{N}$\,: \qquad
$I_rI_s=\left\{
\begin{array}{rcl}
0 &\mbox{if} & r\neq s,\vspace*{2mm} \\
I_r &\mbox{if} & r=s;\\
\end{array}
\right.$

\vskip5mm

2. \,\,  $\forall\,r,s\in[m+1,n]\cap\mathbb{N}$\,: \qquad $I_rI_s=
\sum\limits_{k=\max\{r,s\}+1}^n\Upsilon_{r,k}^{s}I_k$\,;

\vskip5mm

3.\,\, $\forall\,s\in[m+1,n]\cap\mathbb{N}$\;  $\exists!\;
 u_s\in[1,m]\cap\mathbb{N}$ \;$\forall\,
r\in[1,m]\cap\mathbb{N}$\,:\;\;

\begin{equation}\label{mult_rule_3}
I_rI_s=\left\{
\begin{array}{ccl}
0 \;\;\mbox{if}\;\;  r\neq u_s\,,\vspace*{2mm}\\
I_s\;\;\mbox{if}\;\;  r= u_s\,. \\
\end{array}
\right.\medskip
\end{equation}
Moreover, the structure constants
$\Upsilon_{r,k}^{s}\in\mathbb{C}$
 satisfy the associativity conditions:
\vskip2mm (A\,1).\,\, $(I_rI_s)I_p=I_r(I_sI_p)$ \;
$\forall\,r,s,p\in[m+1,n]\cap\mathbb{N}$; \vskip2mm (A\,2).\,\,
$(I_uI_s)I_p=I_u(I_sI_p)$ \; $\forall\,u\in[1,m]\cap\mathbb{N}$\;
 $\forall\,s,p\in[m+1,n]\cap\mathbb{N}$.
\vskip2mm

Obviously, the first $m$ basic vectors $\{I_u\}_{u=1}^m$
are idempotents and 
form a semi-simple subalgebra of the algebra $\mathbb{A}_n^m$. The
vectors $\{I_r\}_{r=m+1}^n$ form a nilpotent subalgebra of the
algebra $\mathbb{A}_n^m$. The element $1=\sum_{u=1}^mI_u$ is the
unit of $\mathbb{A}_n^m$.

In the cases where $\mathbb{A}_n^m$ has some specific properties, 
the following
propositions are true.

\textbf{Proposition 1.} \textit{If there exists the unique
$u_0\in[1,m]\cap\mathbb{N}$ such that $I_{u_0}I_s=I_s$ for all
$s=m+1,\ldots,n$,  then the associativity condition \em (A\,2) \em
is satisfied.}

\textbf{Proof.} In the associativity condition (A\,2), two cases
are possible:
$$
1)\qquad I_u\neq I_{u_0}\,,\;\; \mbox{then}\;\;I_uI_s=0\,\quad
\forall\,s=m+1,\ldots,n;
$$
$$
2)\qquad I_u= I_{u_0}\,,\;\; \mbox{then}\;\;
 I_uI_s=I_s\,\quad \forall\,s=m+1,\ldots,n.
$$

In the first case, the condition
 (A\,2) takes the form
$$0\cdot I_p=I_u\sum\limits_{k=\max\{s,p\}+1}^n\Upsilon_{s,k}^{p}\,I_k=0,
$$
because $I_uI_k=0$ for all $k=\max\{s,p\}+1,\ldots,n$.

In the second case, the condition
 (A\,2) takes the form 
$$I_sI_p=I_u\sum\limits_{k=\max\{s,p\}+1}^n\Upsilon_{s,k}^{p}\,I_k\,.
$$
It is equivalent to the equality 
$$\sum\limits_{k=\max\{s,p\}+1}^n\Upsilon_{s,k}^{p}\,I_k=
\sum\limits_{k=\max\{s,p\}+1}^n\Upsilon_{s,k}^{p}\,I_k,
$$
because $I_uI_k=I_k$ for all $k=\max\{s,p\}+1,\ldots,n$. The
proposition is proved.

Thus, under the conditions of Proposition 1,
the associativity condition (A\,1) is only required. 
It means that the nilpotent subalgebra of $\mathbb{A}_n^m$ with
the basis $\{I_r\}_{r=m+1}^n$  can be an arbitrary commutative
associative nilpotent algebra of dimension $n-m$. Note that such
nilpotent algebras are fully described for the dimensions
$1,2,3,4$ in the paper \cite{Burde_de_Graaf}.

\textbf{Proposition 2.} \textit{If all $u_r$ are different in the
multiplication rule \em 3\em ,
 then $I_sI_p=0$ for all $s,p=m+1,\ldots, n$.}

\textbf{Proof.} Let $s\in[m+1,n]\cap\mathbb{N}$.
 We choose $I_u$ such that $I_uI_s=I_s$.
Then from the associativity condition (A\,2) we obtain the identity
$$I_sI_p=I_u\sum\limits_{k=\max\{s,p\}+1}^n\Upsilon_{s,k}^{p}\,I_k=0,
$$
because by assumption of theorem $I_uI_k=0$ for all
$k=\max\{s,p\}+1,\ldots,n$. The proposition is proved.

Thus, under the conditions of Proposition 2,
the multiplication table of the
nilpotent subalgebra of $\mathbb{A}_n^m$ with the basis
$\{I_r\}_{r=m+1}^n$ consists only of zeros, and all associativity
conditions are satisfied.

The algebra $\mathbb{A}_n^m$ contains $m$ maximal ideals
$$\mathcal{I}_u:=\Biggr\{\sum\limits_{k=1,\,k\neq u}^n\lambda_kI_k:\lambda_k\in
\mathbb{C}\Biggr\}, \quad  u=1,2,\ldots,m,
$$
and their intersection is the radical $$\mathcal{R}:=
\Bigr\{\sum\limits_{k=m+1}^n\lambda_kI_k:\lambda_k\in
\mathbb{C}\Bigr\}.$$

Consider $m$ linear functionals
$f_u:\mathbb{A}_n^m\rightarrow\mathbb{C}$ satisfying the
equalities
$$f_u(I_u)=1,\quad f_u(\omega)=0\quad\forall\,\omega\in\mathcal{I}_u\,,
\quad u=1,2,\ldots,m.
$$
Inasmuch as the kernel of functional $f_u$ is the maximal ideal
$\mathcal{I}_u$, this functional is also continuous and
multiplicative (see \cite[p. 147]{Hil_Filips}).

\vskip3mm

\section{Monogenic functions.}

We consider the vectors $e_1=1,e_2,e_3$ in $\mathbb{A}_n^m$ which
are linearly independent over the field of real numbers
$\mathbb{R}$. It means that the equality
$$\alpha_1e_1+\alpha_2e_2+\alpha_3e_3=0,\quad \alpha_1,\alpha_2,
\alpha_3\in\mathbb{R},$$
holds if and only if $\alpha_1=\alpha_2=
\alpha_3=0$.

Let the vectors $e_1=1,e_2,e_3$ have the following decompositions with respect to
the basis $\{I_k\}_{k=1}^n$:
\begin{equation}\label{e_1_e_2_e_3}
e_1=\sum\limits_{r=1}^mI_r\,, 
\quad e_2=\sum\limits_{r=1}^na_rI_r\,,\quad e_3=\sum\limits_{r=1}^nb_rI_r\,,
\end{equation}
where $a_r,b_r\in\mathbb{C}$.

Let $\zeta:=xe_1+ye_2+ze_3$, where $x,y,z\in\mathbb{R}$. It is
obvious that
 $\xi_u:=f_u(\zeta)=x+ya_u+zb_u$,\, $u=1,2,\ldots,m$.
 Let
 $E_3:=\{\zeta=xe_1+ye_2+ze_3:\,\, x,y,z\in\mathbb{R}\}$ be the
linear span of vectors $e_1,e_2,e_3$ over the field  
$\mathbb{R}$.

Let $\Omega$ be a domain in ${\mathbb R}^3$. Associate with
$\Omega$ the domain $\Omega_{\zeta}:=\{\zeta=xe_1+ye_2+ze_3 :
(x,y,z)\in \Omega\}$ in $E_3$.

We say that a continuous function
$\Phi:\Omega_{\zeta}\rightarrow\mathbb{A}_n^m$ is \textit{monogenic}
in $\Omega_{\zeta}$ if $\Phi$ is differentiable in the sense of
Gateaux in every point of $\Omega_{\zeta}$, i.~e. if  for every
$\zeta\in\Omega_{\zeta}$ there exists an element
$\Phi'(\zeta)\in\mathbb{A}_n^m$ such that
\begin{equation}\label{monogennaOZNA}\medskip
\lim\limits_{\varepsilon\rightarrow 0+0}
\left(\Phi(\zeta+\varepsilon
h)-\Phi(\zeta)\right)\varepsilon^{-1}= h\Phi'(\zeta)\quad\forall\,
h\in E_{3}.\medskip
\end{equation}
$\Phi'(\zeta)$ is the \textit{Gateaux derivative} of the function
$\Phi$ in the point $\zeta$.

In the scientific literature the denomination of monogenic function is used
else for functions satisfying certain conditions similar to the classical Cauchy –
Riemann conditions (see, e.~g., \cite{Brackx, Ryan}).
Such functions are also called
regular functions
(see \cite{Sadberi}) or hyperholomorphic functions (see, e.~g.,
\cite{Kravchenko-Shapiro, Spros}).

Consider the decomposition of a function
$\Phi:\Omega_{\zeta}\rightarrow\mathbb{A}_n^m$ with respect to the
basis $\{I_k\}_{k=1}^n$:
\begin{equation}\label{rozklad-Phi-v-bazysi}
\Phi(\zeta)=\sum_{k=1}^n U_k(x,y,z)\,I_k\,.
 \end{equation}

In the case where the functions $U_k:\Omega\rightarrow\mathbb{C}$ are
$\mathbb{R}$-differentiable in $\Omega$, i.~e. for every
$(x,y,z)\in\Omega$,
$$U_k(x+\Delta x,y+\Delta y,z+\Delta z)-U_k(x,y,z)=
\frac{\partial U_k}{\partial x}\,\Delta x+
\frac{\partial U_k}{\partial y}\,\Delta y+\frac{\partial
U_k}{\partial z}\,\Delta z+$$ $$+\,o\left(\sqrt{(\Delta
x)^2+(\Delta y)^2+(\Delta z)^2}\,\right), \qquad (\Delta
x)^2+(\Delta y)^2+(\Delta z)^2\to 0\,,$$
the function $\Phi$ is monogenic in the domain $\Omega_{\zeta}$ if
and only if the following Cauchy~-- Riemann conditions are
satisfied in $\Omega_{\zeta}$:
\begin{equation}\label{Umovy_K-R}
\frac{\partial \Phi}{\partial y}=\frac{\partial \Phi}{\partial
x}\,e_{2},\quad \frac{\partial \Phi}{\partial z}=\frac{\partial
\Phi}{\partial x}\,e_{3}.
\end{equation}

Below, it will be shown that all components $U_k$ of the monogenic
function \eqref{rozklad-Phi-v-bazysi} are infinitely
$\mathbb{R}$-differentiable in $\Omega$.


\vskip3mm

\section{An expansion of the resolvent.}
\vskip3mm

\lem\label{lem_1_rezolv_A_n_m} \textit{An expansion of the
resolvent is of the form
$$(te_1-\zeta)^{-1}=\sum\limits_{r=1}^nA_rI_r \qquad \forall\,t\in\mathbb{C}:\,
t\neq \xi_u,\quad u=1,2,\ldots,m,
$$
with the coefficients\, $A_r$\, are determined by the following
recurrence relations:
\begin{equation}\label{lem_1_A_p}
\begin{array}{c}
\displaystyle A_u=\frac{1}{t-\xi_u}\,,\;\;u=1,2,\ldots,m, \quad
 A_{m+1}=\frac{T_{m+1}}{\left(t-\xi_{u_{m+1}}\right)^2}\,,\vspace*{4mm}\\
 \displaystyle
A_{p}=\frac{T_p}{(t-\xi_{u_p})^2}+\frac{1}{t-\xi_{u_p}}
\sum\limits_{r=m+1}^{p-1}A_rB_{r,p}\,,\;
\;\;p=m+2,m+3,\ldots,n,\\
\end{array}
\end{equation}
where
\begin{equation}\label{lem_1_T_p}
T_p:=ya_p+zb_p\,,\;\;\;p=m+1,m+2,\ldots,n,
\end{equation}
\begin{equation}\label{lem_1_B_p}
B_{r,p}:=\sum\limits_{s=m+1}^{p-1}T_s \Upsilon_{r,p}^s\,,\;
\;\;p=m+2,m+3,\ldots,n,
\end{equation}
 and the natural numbers $u_p$ are defined in the
rule \em 3 \em of the multiplication table of algebra
$\mathbb{A}_n^m$.}

\textbf{Proof.} Let's find $t\in\mathbb{C}$ for which there exists
the element $(te_1-\zeta)^{-1}$ in the algebra $\mathbb{A}_n^m$
and let's find the coefficients $A_r$ of decomposition of this
element with respect to the basis $\{I_k\}_{k=1}^n$:
$$(te_1-\zeta)^{-1}=\sum\limits_{r=1}^nA_rI_r\,.$$

Taking into account the decompositions (\ref{e_1_e_2_e_3})
and the multiplication table of algebra $\mathbb{A}_n^m$, we obtain 
$$\sum\limits_{u=1}^mI_u=1=(te_1-\zeta)^{-1}(te_1-\zeta)=$$ $$=\Bigr(\sum\limits_{u=1}^mA_uI_u+
\sum\limits_{r=m+1}^nA_rI_r\Bigr)\Bigr(\sum\limits_{u=1}^m(t-\xi_u)I_u-
\sum\limits_{r=m+1}^n(ya_r+zb_r)I_r\Bigr)=
$$
$$=\sum\limits_{u=1}^mA_u(t-\xi_u)I_u+
\Bigr(A_{m+1}(t-\xi_{u_{m+1}})+A_{u_{m+1}}(-ya_{m+1}-zb_{m+1})\Bigr)I_{m+1}+
$$
$$+\sum\limits_{p=m+2}^n\Bigr(A_{u_p}(-ya_p-zb_p)+A_p(t-\xi_{u_p})
+\sum\limits_{r=m+1}^{p-1}A_r\sum\limits_{s=m+1}^{p-1}(-ya_s-zb_s)\Upsilon_{r,p}^s
\Bigr)I_p\,.
$$


Now, to determine the coefficients $A_r$, we have the system of
equations:
$$
A_u(t-\xi_u)=1,\;\;u=1,2,\ldots,m, \;\;
 A_{m+1}(t-\xi_{u_{m+1}})-A_{u_{m+1}}T_{m+1}=0,$$
 $$-A_{u_p}T_p+A_p(t-\xi_{u_p})
-\sum\limits_{r=m+1}^{p-1}A_rB_{r,p}=0,\;\; p=m+2,m+3,\ldots,n,
$$
where the denotations  (\ref{lem_1_T_p}), (\ref{lem_1_B_p}) are
used.
As immediate consequences of this system, we obtain the equalities
(\ref{lem_1_A_p}). The lemma is proved.

It follows from Lemma \ref{lem_1_rezolv_A_n_m}  that the points
 $(x,y,z)\in\mathbb{R}^3$
corresponding to the noninvertible elements
$\zeta=xe_1+ye_2+ze_3$ 
 form the straight lines
  \[L_u:\quad\left\{
\begin{array}{r}x+y{\rm Re}\,a_u+z{\rm Re}\,b_u=0,\vspace*{3mm} \\
y{\rm Im}\,a_u+z{\rm Im}\,b_u=0, \\ \medskip
\end{array} \right. \qquad u=1,2,\ldots,m,\]
in the three-dimensional space $\mathbb{R}^3$.

\lem\label{lem_2_rezolv_A_n_m} \textit{If there exists $p\in[m+2,n]\cap\mathbb{N}$ and
$r\in[m+1,p-1]\cap\mathbb{N}$ such that
$B_{r,p}\neq 0$,
 then $u_p=u_r$.}

\textbf{Proof.}  Since $B_{r,p}\neq 0$, at least one of the
numbers
 $\Upsilon_{m+1,p}^{r}$, $\Upsilon_{m+2,p}^{r},\ldots,\Upsilon_{p-1,p}^{r}$
is different from zero. Let  $\Upsilon_{k,p}^r\neq0$, where $k$ is
one of the numbers $m+1,m+2,\ldots,p-1$. The associativity
condition of the algebra implies the equality
$$(I_{u_{r}}I_r)I_k=I_{u_{r}}(I_rI_k),
$$
which is equivalent to the equality
\begin{equation}\label{lem_2_1}
\sum\limits_{\ell=\max\{k,r\}+1}^n\Upsilon_{k,\ell}^rI_\ell=
\sum\limits_{\ell=\max\{k,r\}+1}^n\Upsilon_{k,\ell}^rI_{u_{r}}I_\ell\,.
\end{equation}

Note that by the multiplication rule (\ref{mult_rule_3}), for each
$\ell=m+1,m+2,\ldots,n$ the product $I_{u_{r}}I_\ell$ is equal to
either zero or $I_\ell$. Therefore, since
 $\Upsilon_{k,p}^r\neq0$, 
the equality (\ref{lem_2_1}) implies the equality
$\Upsilon_{k,p}^rI_p=\Upsilon_{k,p}^rI_{u_{r}}I_p$\,,
 i.~e. $I_p=I_{u_{r}}I_p$, and it means that $u_r=u_p$. The lemma is proved.


\lem\label{lem_3_rezolv_A_n_m} \textit{For any
$s=m+1,m+2,\ldots,n$ the coefficients
 $A_s$
can be 
expressed in the form
\begin{equation}\label{lem_3_1}
A_s=\sum\limits_{k=2}^{s-m+1}\frac{Q_{k,s}}{\left(t-\xi_{u_s}\right)^k}\,,
\end{equation}
where $Q_{k,s}$ are determined by the following recurrence
relations:}
\begin{equation}\label{lem_3_2}
\begin{array}{c}
\displaystyle
Q_{2,s}:=T_{s}\,,\qquad
Q_{k,s}=\sum\limits_{r=m+1}^{s-1}Q_{k-1,r}\,B_{r,\,s}\,,\; \;\;k=3,4,\ldots,s-m+1.\\
\end{array}
\end{equation}

\textbf{Proof.} Let us prove the expression (\ref{lem_3_1}) by
mathematical induction. First, 
the expression (\ref{lem_3_1}) coincides with the equality
(\ref{lem_1_A_p}) for $s=m+1$.


Further, suppose the correctness of the formula (\ref{lem_3_1})
for all $A_{m+1}$, $A_{m+2},\ldots,A_{s-1}$ and prove that $A_s$\,
is also expressed by the formula (\ref{lem_3_1}). For this
purpose, we use the equality (\ref{lem_1_A_p}) for $p=s$.
 Substituting the expressions (\ref{lem_3_1}) for $A_r$
 in the equality (\ref{lem_1_A_p}), we obtain:\medskip
 $$A_s=\frac{T_s}{\left(t-\xi_{u_s}\right)^2}+ \sum\limits_{r=m+1}^{s-1}
\frac{A_r\,B_{r,\,s}}{t-\xi_{u_s}}=$$
\begin{equation}\label{lem_3_4}
=\frac{T_s}{\left(t-\xi_{u_s}\right)^2}+
\sum\limits_{r=m+1}^{s-1}\sum\limits_{k=2}^{r-m+1}\frac{Q_{k,r}\,B_{r,\,s}}
{(t-\xi_{u_s})(t-\xi_{u_r})^k}\,.
\end{equation}\vskip4mm
If all $B_{r,\,s}=0$ for $r=m+1,m+2,\ldots,s-1$, then the formula
(\ref{lem_3_4}) takes the form (\ref{lem_3_1}) with
$Q_{2,m+2}=T_{m+2}$ and $Q_{k,s}=0$. Furthermore, for every
$r=m+1,m+2,\ldots,s-1$ for which $B_{r,\,s}\neq0$, by Lemma
\ref{lem_2_rezolv_A_n_m}\, $u_r=u_s$, and we get again the formula
(\ref{lem_3_1}), where $Q_{k,s}$ are determined by the equalities
(\ref{lem_3_2}). The lemma is proved.

As a consequence of Lemmas \ref{lem_1_rezolv_A_n_m},
\ref{lem_3_rezolv_A_n_m}, we obtain the following expansion of
resolvent:
 \begin{equation}\label{lem_3_5}
 (te_1-\zeta)^{-1}=\sum\limits_{u=1}^m\frac{1}{t-\xi_u}\,I_u+
 \sum\limits_{s=m+1}^{n}\sum\limits_{k=2}^{s-m+1}\frac{Q_{k,s}}
 {\left(t-\xi_{u_{s}}\right)^k}\,I_{s}\,.
 \end{equation}
\vskip3mm

\section{A constructive description of monogenic functions.}
\vskip3mm

Denote $f_u(E_3):=\{f_u(\zeta) : \zeta\in E_3\}$. In what follows,
we make the following essential assumption:
$f_u(E_3)=\mathbb{C}$ for all\, $u=1,2,\ldots,m$.
 Obviously, it holds if and only if for every fixed $u=1,2, \ldots, m$
at least one of the numbers $a_u$ or $b_u$ belongs to
$\mathbb{C}\setminus\mathbb{R}$.

\lem\label{lem_1_konstruct_opys_A_n_m} \textit{Suppose that a
domain $\Omega\subset \mathbb{R}^{3}$ is convex in the direction
of the straight lines $L_u$ and $f_u(E_3)=\mathbb{C}$ for all
 $u=1,2,\ldots, m$. Suppose also that a
function $\Phi:\Omega_{\zeta}\rightarrow\mathbb{A}_n^m$ is
monogenic in the domain
  $\Omega_{\zeta}$. If points
$\zeta_{1},\zeta_{2}\in\Omega_{\zeta}$ such that
$\zeta_{2}-\zeta_{1}\in \{\zeta=xe_1+ye_2+ze_3:(x,y,z)\in L_u\}$,
then}
\begin{equation}\label{Fi(dz1)-Fi(dz2)}
\Phi(\zeta_{1})-\Phi(\zeta_{2})\in\mathcal{I}_u.
\end{equation}

The proof 
of Lemma \ref{lem_1_konstruct_opys_A_n_m} is similar to the proof
of Lemma 1 \cite{Pl-Shp1}, where one 
must take the straight line $L_u$ and the functional $f_u$ instead
of the straight line $L$ and the functional $f$, respectively.

Let a domain $\Omega\subset \mathbb{R}^{3}$ be convex in the
direction of the straight lines $L_u$,\, $u=1,2,\ldots, m$. By
$D_u$ we denote that domain in $\mathbb{C}$
 onto which the domain $\Omega_\zeta$ is mapped by the functional  $f_u$.

We introduce the linear operators $A_u$,\,$u=1,2,\ldots,m$, which
assign holomorphic functions $F_u:\,D_u\rightarrow\mathbb{C}$ to
every monogenic function
$\Phi:\Omega_{\zeta}\rightarrow\mathbb{A}_n^m$ by the formula
\begin{equation}\label{def_op_A}
F_u(\xi_u)=f_u(\Phi(\zeta)),
\end{equation}
where $\xi_u=f_u(\zeta)\equiv x+ya_u+zb_u$ and
$\zeta\in\Omega_{\zeta}$. It follows from Lemma
\ref{lem_1_konstruct_opys_A_n_m}  that the value $F_u(\xi_u)$ does
not depend on a choice of a point $\zeta$ for which
$f_u(\zeta)=\xi_u$.

\vskip3mm \lem\label{lem_2_konstruct_opys_A_n_m} \textit{Suppose
that a domain $\Omega\subset \mathbb{R}^{3}$ is convex in the
direction of the straight lines $L_u$ and $f_u(E_3)=\mathbb{C}$
for all $u=1,2,\ldots, m$. Suppose also that for any fixed
$u=1,2,\ldots,m$, a function $F_u:D_u\rightarrow \mathbb{C}$ is
holomorphic in a domain $D_u$ and $\Gamma_u$ is a closed Jordan
rectifiable curve in $D_u$ which surrounds the point $\xi_u$ and
contains no points $\xi_q$, $q=1,2,\ldots, m$,\, $q\neq u$. Then
the function
\begin{equation}\label{lem_5}
\Psi_u(\zeta):=I_u\int\limits_{\Gamma_u}F_u(t)(te_1-\zeta)^{-1}\,dt
\end{equation}
is monogenic in the domain
 $\Omega_\zeta$. }

\textbf{Proof.} Let $\zeta\in\Omega_\zeta$. First, for any $h\in
E_3$ and any $\varepsilon>0$,
it follows from the equality (\ref{lem_3_5}) that 
$$
 I_u(te_1-\zeta-\varepsilon h)^{-1}=\sum\limits_{u=1}^m\frac{1}{t-\xi_u-\varepsilon
 f_{u}(h)}\,I_u+$$
 \begin{equation}\label{lem_5_1}
 +\sum\limits_{s=m+1}^{n}\sum\limits_{k=2}^{s-m+1}\frac{Q_{k,s}}
 {\left(t-\xi_{u_s}-\varepsilon
 f_{u_{s}}(h)\right)^k}\,I_{s}\,I_u\,. 
 \end{equation}

Note that for any 
 natural $k$, a complex-value function
 $\displaystyle 1/\left(t-\xi_u-\varepsilon f_u(h)\right)^k$ tends to the function $\displaystyle 1/(t-\xi_u)^{k}$
 uniformly for all $t\in\Gamma_u$ when $\varepsilon\rightarrow0$.
Therefore, for any $h\in E_3$ the function
$I_u(te_1-\zeta-\varepsilon h)^{-1}$ tends to the function
$I_u(te_1-\zeta)^{-1}$ uniformly for all $t\in\Gamma_u$ when
$\varepsilon\rightarrow0$.

Further, let us prove the existence of Gateaux derivative 
$\Psi'_u(\zeta)$ by the definition (\ref{monogennaOZNA}).
Taking into account the Hilbert identity (see, e.~g., Theorem
4.8.2 \cite[p.\,140]{Hil_Filips})
$$(te_1-a)^{-1}-(te_1-b)^{-1}=(te_1-a)^{-1}(te_1-b)^{-1}(a-b)\qquad \forall\,a,b\in E_3\,,
$$ we have
$$\Lambda:=I_u\lim\limits_{\varepsilon\rightarrow 0+0}\frac{1}{\varepsilon}
\Biggr(\int\limits_{\Gamma_u}F_u(t)\big(te_1-(\zeta+\varepsilon
h)\big)^{-1}\,dt-
\int\limits_{\Gamma_u}F_u(t)(te_1-\zeta)^{-1}\,dt\Biggr)=
 $$
$$=I_uh\lim\limits_{\varepsilon\rightarrow 0+0}\int\limits_{\Gamma_u}F_u(t)
\big(te_1-(\zeta+\varepsilon h)\big)^{-1}(te_1-\zeta)^{-1}\,dt.
$$
Now, since the function $I_u(te_1-\zeta-\varepsilon h)^{-1}$
converges uniformly to the function $I_u(te_1-\zeta)^{-1}$ when
$\varepsilon\rightarrow0$, we obtain the equality
$$\Lambda=I_uh\int\limits_{\Gamma_u}F_u(t)\big((te_1-\zeta)^{-1}\big)^2\,dt 
$$
which means the existence of the Gateaux derivative\,
$$\Psi_u'(\zeta):=I_u\int\limits_{\Gamma_u}F_u(t)\big((te_1-\zeta)^{-1}\big)^2\,dt.$$

Finally, in view of the equality (\ref{lem_3_5}),
components of the expansion of function $\Psi_u(\zeta)$ with
respect to the basis $\{I_k\}_{k=1}^{n}$ are continuous functions.
Therefore, the function $\Psi'_u(\zeta)$ is also continuous, and
the function $\Psi'_u(\zeta)$ is monogenic in $\Omega_\zeta$.
Lemma is proved.

\lem\label{Lem-6-osn-const-op} \textit{Suppose that a domain
$\Omega\subset\mathbb{R}^{3}$ is convex in the direction of the
straight lines $L_u$ and $f_u(E_3)=\mathbb{C}$ for all
 $u=1,2,\ldots, m$. Suppose also that a function $V:\Omega\rightarrow\mathbb{C}$
 satisfies the equalities
 \begin{equation}\label{lem-6-1}
\frac{\partial V}{\partial y}=\frac{\partial V}{\partial x}\,
a_{u}\,, \qquad \frac{\partial V}{\partial z}=\frac{\partial
V}{\partial x}\, b_{u}
\end{equation}
in $\Omega$. Then $V$ is a holomorphic function of the variable
$\xi_u=f_u(\zeta)=x+ya_u+zb_u$ in the domain $D_u$.}

\textbf{Proof.} We first separate the real and the imaginary part
of the expression
\begin{equation}\label{lem-6-2}
\xi_u=x+y\,{\rm Re}\,a_u+z\,{\rm Re}\,b_u + i(y\,{\rm
Im}\,a_u+z\,{\rm Im}\,b_u)=: \tau+ i\eta
\end{equation}
and note that the equalities (\ref{lem-6-1}) yield
\begin{equation}\label{lem-6-3}
\frac{\partial V}{\partial \eta}\,{\rm Im}\,a_u =i\,\frac{\partial
V}{\partial \tau}\,{\rm Im}\,a_u\,, \qquad \frac{\partial
V}{\partial \eta}\,{\rm Im}\,b_u=i\,\frac{\partial V}{\partial
\tau}\, {\rm Im}\,b_u.
\end{equation}

It follows from the condition $f_u(E_3)=\mathbb{C}$ that
 at least one of the numbers ${\rm Im}\,a_u$ or ${\rm Im}\,b_u$ is not
equal to zero. Therefore, using (\ref{lem-6-3}), we get
\begin{equation}\label{lem-6-4}
\frac{\partial V}{\partial \eta}=i\,\frac{\partial V}{\partial
\tau}\,.
\end{equation}

We prove that $V(x_1, y_1, z_1)=V(x_2, y_2, z_2)$ for points
$(x_1,y_1,z_1),(x_2,y_2,z_2)\in\Omega$ such that the segment that
connects these points is parallel to the straight line $L_u$. To
this end, in the domain $\Omega$, we construct two surfaces with
common edge,
 namely a surface $Q$ that contains the point $(x_1,y_1,z_1)$ and a surface
 $\Sigma$  that contains the point $(x_2,y_2,z_2)$, such that the
  restrictions of the functional
 $f_u$ to the corresponding subsets $Q_\zeta:=\{\zeta\in E_3:(x,y,z)\in Q\}$ and
 $\Sigma_\zeta:=\{\zeta\in E_3:(x,y,z)\in \Sigma\}$ of the domain $\Omega_\zeta$
  are bijections of these subsets to the same domain $D_u$ of the complex plane.

As the surface $Q$ in the domain $\Omega$, we take a fixed
equilateral triangle
 with vertices $A_{u,1},A_{u,2}$, and $A_{u,3}$
centered at the point $(x_1,y_1,z_1)$ the plane of which is
perpendicular to the straight line $L_u$. We now continue the
construction of the surface $\Sigma$.

Consider the triangle with vertices $A_{u,1}'$, $A_{u,2}'$, and
$A_{u,3}'$
 centered at the
point $(x_2,y_2,z_2)$, lying in the domain $\Omega$, and such that
its sides $A_{u,1}'A_{u,2}'$, $A_{u,2}'A_{u,3}'$, and
$A_{u,1}'A_{u,3}'$ are parallel to the segments $A_{u,1}A_{u,2}$,
$A_{u,2}A_{u,3}$, and $A_{u,1}A_{u,3}$,
 respectively,
and have smaller lengths than the sides of the triangle
$A_{u,1}A_{u,2}A_{u,3}$. Since the domain $\Omega$ is convex in
the direction of the straight line $L_u$, we conclude that the
prism with vertices $A_{u,1}'$, $A_{u,2}'$, $A_{u,3}'$,
 $A_{u,1}''$, $A_{u,2}''$,
 and $A_{u,3}''$ such that the points $A_{u,1}''$, $A_{u,2}''$, and $A_{u,3}''$
  lie in the plane
 of the triangle $A_{u,1}A_{u,2}A_{u,3}$ and its edges $A_{u,m}'A_{u,m}''$,
 $m=\overline{1,3}$,
  are parallel to the straight line $L_u$ is completely contained in $\Omega$.

We now fix a triangle with vertices $B_{u,1}$, $B_{u,2}$, and
$B_{u,3}$
 such that the point
$B_{u,m}$ lies on the segment $A_{u,m}'A_{u,m}''$ for
$m=\overline{1,3}$
 and the truncated pyramid with vertices $A_{u,1}$, $A_{u,2}$, $A_{u,3}$, $B_{u,1}$,
  $B_{u,2}$,  and $B_{u,3}$ and lateral edges $A_{u,m}B_{u,m}$,
$m=\overline{1,3}$, is completely contained in the domain
$\Omega$.

Finally, in the plane of the triangle $A_{u,1}'A_{u,2}'A_{u,3}'$,
we fix a triangle $T$  with vertices $C_{u,1}$, $C_{u,2}$, and
$C_{u,3}$
 such that its sides $C_{u,1}C_{u,2}$, $C_{u,2}C_{u,3}$,
 and $C_{u,1}C_{u,3}$ are parallel to the segments $A_{u,1}'A_{u,2}'$,
  $A_{u,2}'A_{u,3}'$, and $A_{u,1}'A_{u,3}'$,
  respectively, and have smaller lengths than the sides of the triangle
 $A_{u,1}'A_{u,2}'A_{u,3}'$. By construction, the truncated pyramid with
vertices $B_{u,1}$, $B_{u,2}$, $B_{u,3}$, $C_{u,1}$, $C_{u,2}$,
and $C_{u,3}$
 and lateral edges $B_{u,m}C_{u,m}$,
 $m=\overline{1,3}$, is completely contained in the domain $\Omega$.

Let $\Sigma$ denote the surface formed by the triangle $T$ and the
lateral
 surfaces of the truncated pyramids $A_{u,1}A_{u,2}A_{u,3}B_{u,1}B_{u,2}B_{u,3}$ and $B_{u,1}B_{u,2}B_{u,3}C_{u,1}C_{u,2}C_{u,3}$.

Since the surfaces $Q$ and $\Sigma$ have a common edge, the sets
$Q_\zeta$ and $\Sigma_\zeta$ are mapped by the functional $f_u$
onto the same domain $D_u$ of the complex plane. In the domain
$D_u$, we introduce two complex-valued functions $H_u$ and $W_u$
in $D_u$ as follows:
$$H_u(\xi_u)=V(x,y,z)\quad \text{for}\quad (x,y,z)\in Q,$$
$$W_u(\xi_u)=V(x,y,z)\quad \text{for} \quad(x,y,z)\in\Sigma,$$
where the correspondence between the points $(x,y,z)$ and
$\xi_u\in D_u$ is described by relation (\ref{lem-6-2}).

By virtue of equality (\ref{lem-6-4}) and Theorem 6 in
\cite{Tolstov},
 the functions $H_u$ and $W_u$ are holomorphic in the domain $D_u$.
 According to the definition of the functions $H_u$ and $W_u$,
  we have $H_u(\xi_u)\equiv W_u(\xi_u)$ on the boundary of
the domain $D_u$. By virtue of the holomorphy of the functions
$H_u$ and $W_u$ in the domain $D_u$, the identity
$H_u(\xi_u)\equiv W_u(\xi_u)$ holds everywhere in $D_u$.
Therefore, the equality $V(x_1, y_1, z_1)=V(x_2, y_2, z_2)$ is
proved.

Thus, a function $V:\Omega\rightarrow\mathbb{C}$ of the form
$V(x,y,z):=F(\xi_u)$,
 where $F(\xi_u)$ is an arbitrary function holomorphic in
the domain $D_u$, is a general solution of the system
(\ref{lem-6-1}). The lemma is proved.

\vskip3mm \theor\label{teo_1_konstruct_opys_A_n_m} \textit{Let a
domain $\Omega\subset \mathbb{R}^{3}$ be convex in the direction
of the straight lines $L_u$ and $f_u(E_3)=\mathbb{C}$ for all
 $u=1,2,\ldots, m$.  Then every
monogenic function $\Phi:\Omega_{\zeta}\rightarrow\mathbb{A}_n^m$
can be expressed in the form
 \begin{equation}\label{Teor--1}
\Phi(\zeta)=\sum\limits_{u=1}^mI_u\,\frac{1}{2\pi
i}\int\limits_{\Gamma_u} F_u(t)(te_1-\zeta)^{-1}\,dt+
\sum\limits_{s=m+1}^nI_s\,\frac{1}{2\pi i}\int\limits_
{\Gamma_{u_s}}G_s(t)(te_1-\zeta)^{-1}\,dt,
 \end{equation}
where $F_u$ and $G_s$ are certain holomorphic functions in the
domains $D_u$ and $D_{u_s}$, respectively, and $\Gamma_q$ is a
closed Jordan rectifiable curve in $D_q$ which surrounds the point
$\xi_q$ and contains no points $\xi_{\ell}$, $\ell,q=1,2,\ldots,
m$,\,$\ell\neq q$.}

\textbf{Proof.} We set
 \begin{equation}\label{teor__1}
 F_u:=A_u\Phi,\;\;\;u=1,2,\ldots,m.
 \end{equation}
Let us show that the values of monogenic function
\begin{equation}\label{teor__2}
\Phi_0(\zeta):=\Phi(\zeta)-\sum\limits_{u=1}^mI_u\,\frac{1}{2\pi
i} \int\limits_{\Gamma_u}F_u(t)(te_1-\zeta)^{-1}\,dt
\end{equation}
belong to the radical $\mathcal{R}$, 
i.~e. $\Phi_0(\zeta)\in\mathcal{R}$ for all
$\zeta\in\Omega_\zeta$. 
As a consequence of the equality (\ref{lem_3_5}), we have the
equality
$$I_u\,\frac{1}{2\pi i}
\int\limits_{\Gamma_u}F_u(t)(te_1-\zeta)^{-1}\,dt=I_u\,\frac{1}{2\pi
i} \int\limits_{\Gamma_u}\frac{F_u(t)}{t-\xi_u}\,dt+$$
$$
 +\frac{1}{2\pi i}\sum\limits_{s=m+1}^{n}
 \sum\limits_{k=2}^{s-m+1}\int\limits_{\Gamma_u}\frac{F_u(t)Q_{k,s}}
 {\left(t-\xi_{u_{s}}\right)^k}\,dt
 \,I_{s}\,I_u\,,
$$
from which we obtain the equality
\begin{equation}\label{teor__3}
f_u\left(\sum\limits_{u=1}^mI_u\,\frac{1}{2\pi i}
\int\limits_{\Gamma_u}F_u(t)(te_1-\zeta)^{-1}\,dt\right)=F_u(\xi_u).
\end{equation}
Operating onto the  equality (\ref{teor__2}) by the functional
 $f_u$ and taking into account the relations (\ref{def_op_A}), (\ref{teor__1}), (\ref{teor__3}),  we get the equality
$$f_u(\Phi_0(\zeta))=F_u(\xi_u)-F_u(\xi_u)=0
$$
for all $u=1,2,\ldots,m$, i.~e. $\Phi_0(\zeta)\in\mathcal{R}$.

Therefore, the function $\Phi_0$ is of the form
\begin{equation}\label{Fi_0--1}
\Phi_{0}(\zeta)=\sum\limits_{s=m+1}^{n} V_{s}(x,y,z)\,I_s\,,
\end{equation}
where $V_{s}:\Omega\rightarrow\mathbb{C}$\,, and the Cauchy~--
Riemann conditions (\ref{Umovy_K-R}) are satisfied with
$\Phi=\Phi_0$.
Substituting the expressions (\ref{e_1_e_2_e_3}), (\ref{Fi_0--1}) into the equality
(\ref{Umovy_K-R}), we obtain
\begin{equation}\label{teor__5-1}
\begin{array}{c}
\displaystyle
\sum\limits_{s=m+1}^{n} \frac{\partial V_{s}}{\partial y}\,I_s=
\sum\limits_{s=m+1}^{n} \frac{\partial V_{s}}{\partial x}\,I_s
\sum\limits_{r=1}^n a_r\,I_r\,,
\vspace*{4mm}\\
\displaystyle
\sum\limits_{s=m+1}^{n} \frac{\partial V_{s}}{\partial z}\,I_s=
\sum\limits_{s=m+1}^{n} \frac{\partial V_{s}}{\partial x}\,I_s
\sum\limits_{r=1}^n b_r\,I_r\,.\\
\end{array}
\end{equation}
Equating the coefficients of $I_{m+1}$ in these equalities,
 we obtain the following system of equations for determining the function
 $V_{m+1}(x,y,z)$:
$$\frac{\partial V_{m+1}}{\partial y}=\frac{\partial V_{m+1}}{\partial x}\,
a_{u_{m+1}}\,,
\quad
\frac{\partial V_{m+1}}{\partial z}=\frac{\partial V_{m+1}}{\partial x}\,
b_{u_{m+1}}\,.
$$
It follows from Lemma \ref{Lem-6-osn-const-op} that
$V_{m+1}(x,y,z)\equiv G_{m+1}(\xi_{u_{m+1}})$, where $G_{m+1}$ is
a function holomorphic in the domain $D_{u_{m+1}}$. Therefore,
\begin{equation}\label{teor__4}
\Phi_0(\zeta)=G_{m+1}(\xi_{u_{m+1}})\,I_{m+1}+
\sum\limits_{s=m+2}^{n} V_{s}(x,y,z)\,I_s.
\end{equation}

Due to the expansion (\ref{lem_3_5}), we have the representation
\begin{equation}\label{teor__5}
I_{m+1}\,\frac{1}{2\pi
i}\int\limits_{\Gamma_{u_{m+1}}}G_{m+1}(t)(te_1-\zeta)^{-1}\,dt=
G_{m+1}(\xi_{u_{m+1}})\,I_{m+1}+\Psi(\zeta),
\end{equation}
where $\Psi(\zeta)$ is a function with values in the set
$\big\{\sum_{k=m+2}^n\alpha_k\,I_k:\alpha_k\in\mathbb{C}\big\}$.

Now, consider the function
 $$\Phi_1(\zeta):=\Phi_0(\zeta)-I_{m+1}\,\frac{1}{2\pi
i}\int\limits_{\Gamma_{u_{m+1}}}G_{m+1}(t)(te_1-\zeta)^{-1}\,dt.$$
In view of the relations (\ref{teor__4}), (\ref{teor__5}),
$\Phi_1$ can be represented in the form
$$\Phi_{1}(\zeta)=\sum\limits_{s=m+2}^n \widetilde
V_{s}(x,y,z)\,I_s,$$
where $\widetilde V_{s}:\Omega\rightarrow\mathbb{C}$\,.

Inasmuch as $\Phi_1$ is a monogenic function in $\Omega_{\zeta}$,
the functions $\widetilde V_{m+2},\widetilde
V_{m+3},\dots,\widetilde V_{n}$ satisfy the system
\eqref{teor__5-1}, where $V_{m+1}\equiv 0$, $V_s=\widetilde V_{s}$
for $s=m+2,m+3,\ldots,n$. Therefore,  similarly to the function
$V_{m+1}(x,y,z)\equiv G_{m+1}(\xi_{u_{m+1}})$, the function
$\widetilde V_{m+2}$ satisfies the equations
$$
\frac{\partial \widetilde V_{m+2}}{\partial y}=\frac{\partial \widetilde V_{m+2}}
{\partial x}\,a_{u_{m+2}}\,,\qquad
\frac{\partial \widetilde V_{m+2}}{\partial z}=\frac{\partial
\widetilde V_{m+2}}{\partial x}\,b_{u_{m+2}}
$$
 and is of the
form $\widetilde V_{m+2}(x,y,z)\equiv G_{m+2}(\xi_{u_{m+2}})$,
where $G_{m+2}$ is a function holomorphic in the domain
$D_{u_{m+2}}$.

In such a way, step by step, considering the functions
$$\Phi_j(\zeta):=\Phi_{j-1}(\zeta)-I_{m+j}\,\frac{1}{2\pi
i}\int\limits_{\Gamma_{u_{m+j}}}G_{m+j}(t)(te_1-\zeta)^{-1}\,dt$$
 for $j=2,3,\dots,n-m-1$, we get the representation
(\ref{Teor--1}) of the function $\Phi$. The theorem is proved.

Taking into account the expansion (\ref{lem_3_5}), one can rewrite
the equality (\ref{Teor--1}) in the following equivalent form:
$$\Phi(\zeta)=\sum\limits_{u=1}^mF_u(\xi_u)I_u+
\sum\limits_{s=m+1}^{n}\sum\limits_{k=2}^{s-m+1}\frac{1}{(k-1)!}
Q_{k,s}F_{u_s}^{(k-1)}(\xi_{u_s})\,I_{s}+$$
\begin{equation}\label{dopolnenije-1-1-9}+\sum\limits_{q=m+1}^nG_q(\xi_{u_q})I_q+
\sum\limits_{q=m+1}^n\sum\limits_{s=m+1}^{n}
\sum\limits_{k=2}^{s-m+1}\frac{1}{(k-1)!}Q_{k,s}G_q^{(k-1)}(\xi_{u_q})\,I_{q}
\,I_s\,.
\end{equation}\vskip4mm

Thus, the equalities (\ref{Teor--1}) and (\ref{dopolnenije-1-1-9})
specify methods to construct explicitly any monogenic functions
$\Phi:\Omega_\zeta\rightarrow \mathbb{A}_n^m$ using $n$
corresponding holomorphic functions of complex variables.

The following statement follows immediately from the equality
(\ref{dopolnenije-1-1-9}) in which the right-hand side is a
monogenic function in the domain $\Pi_\zeta:=\{\zeta\in
E_3:f_u(\zeta)=D_u,\,u=1,2,\ldots,m\}$.

\theor\label{teo_pro_naslidky} {\it Let a domain $\Omega$ be
convex in the directions of the straight lines $L_u$ and
$f_u(E_3)=\mathbb{C}$ for all
 $u=1,2,\ldots,m$. Then every monogenic function
$\Phi:\Omega_{\zeta}\rightarrow \mathbb{A}_n^m$ can be continued
to a function monogenic in the domain $\Pi_{\zeta}$.}

The next statement is a fundamental consequence of the equality
(\ref{dopolnenije-1-1-9}), and it is true for an arbitrary domain
$\Omega_{\zeta}$.

\theor\label{teo_pro_naslidky2} {\it Let $f_u(E_3)=\mathbb{C}$ for all
 $u=1,2,\ldots,m$. Then for every monogenic function
$\Phi:\Omega_{\zeta}\rightarrow \mathbb{A}_n^m$ in an arbitrary
domain $\Omega_{\zeta}$, the Gateaux $r$-th derivatives
$\Phi^{(r)}$ are monogenic functions in $\Omega_{\zeta}$ for all\,
$r$.}

The proof is completely analogous to the proof of Theorem 4
\cite{Pl-Shp1}.

 Using the integral expression
(\ref{Teor--1}) of monogenic
function $\Phi:\Omega_{\zeta}\rightarrow \mathbb{A}_n^m$ in the case
where a domain $\Omega$ is convex in the
directions of the straight lines $L_u$\,,\;$u=1,2,\ldots,m$, we obtain the following
expression for the Gateaux $r$-th derivative $\Phi^{(r)}$:
$$\Phi^{(r)}(\zeta)=\sum\limits_{u=1}^mI_u\,\frac{r!}{2\pi i}\int\limits_{\Gamma_u}
F_u(t)\Big((t-\zeta)^{-1}\Big)^{r+1}\,dt+$$
$$
+\sum\limits_{s=m+1}^nI_s\,\frac{r!}{2\pi i}\int\limits_
{\Gamma_{u_s}}G_s(t)\Big((t-\zeta)^{-1}\Big)^{r+1}\,dt\qquad
\forall\zeta\in\Omega_{\zeta}\,.\medskip
$$

\vskip3mm
\section{Special cases.}
\vskip3mm

Note that in the cases where the algebra $\mathbb{A}_n^m$ has some
specific properties (for instance, properties described in
Propositions 1 and 2),
it is easy to simplify the form of the equality
(\ref{dopolnenije-1-1-9}).

\textbf{1.} In the case considered in Proposition 1, 
the following equalities hold:
$$u_{m+1}=u_{m+2}=\ldots=u_n=:\eta\,.$$

In this case the representation (\ref{dopolnenije-1-1-9}) takes
the form
$$\Phi(\zeta)=\sum\limits_{u=1}^mF_u(\xi_u)I_u+
\sum\limits_{s=m+1}^{n}\sum\limits_{k=2}^{s-m+1}\frac{1}{(k-1)!}
Q_{k,s}F_\eta^{(k-1)}(\xi_\eta)\,I_{s}+$$
\begin{equation}\label{dopolnenije-1-1}+\sum\limits_{s=m+1}^nG_s(\xi_\eta)I_s+
\sum\limits_{q=m+1}^n\sum\limits_{s=m+1}^{n}
\sum\limits_{k=2}^{s-m+1}\frac{1}{(k-1)!}Q_{k,s}G_q^{(k-1)}(\xi_\eta)\,I_{s}
\,I_q\,.
\end{equation}\vskip4mm

The formula (\ref{dopolnenije-1-1}) generalizes representations of
monogenic functions in both three-dimensional harmonic algebras
(see \cite{Pl-Shp1,Pl-Pukh,Pukh}) and specific $n$-dimensional
algebras (see \cite{Pl-Shp-Algeria,Pl-Pukh-Analele}) to the case
of algebras more general form.

\vskip3mm \textbf{2.} In the case considered in Proposition 2, 
the functions $B_{r, p}$ 
from the equalities (\ref{lem_1_B_p}) are identically equal to
zero.
 In this case the representation (\ref{lem_3_5}) takes the form
\begin{equation}\label{dopolnenije-2}
 (te_1-\zeta)^{-1}=\sum\limits_{u=1}^m\frac{1}{t-\xi_u}\,I_u+
 \sum\limits_{s=m+1}^{n}\frac{T_{s}}{\left(t-\xi_{u_{s}}\right)^2}\,I_{s}\,,
 \end{equation}
and as a consequence of the equalities (\ref{Teor--1}),
(\ref{dopolnenije-2}), we obtain the following representation of
monogenic function:
 \begin{equation}\label{dopolnenije-3}
\Phi(\zeta)=\sum\limits_{u=1}^mF_u(\xi_u)I_u+\sum\limits_{s=m+1}^nG_s(\xi_{u_s})I_s+
\sum\limits_{s=m+1}^nT_sF_{u_s}^{\,'}(\xi_{u_s})I_s\,.
\end{equation}

The formula (\ref{dopolnenije-3}) generalizes representations of
monogenic functions in both a three-dimensional harmonic algebra
with one-dimensional radical (see \cite{Pl-Pukh}) and semi-simple
algebras (see \cite{Pukh,Pl-Pukh-Analele}) to the case of algebras
more general form.

\vskip3mm \textbf{3.} In the case where $n=m$, the algebra
$\mathbb{A}_n^n$ is semi-simple and contains no nilpotent
subalgebra. Then the formulas (\ref{dopolnenije-1-1}),
(\ref{dopolnenije-3}) take the form
 $$
 \Phi(\zeta)=\sum\limits_{u=1}^nF_u(\xi_u)I_u\,,
$$
because there are no vectors $\{I_k\}_{k=m+1}^n$.
This formula was obtained in the paper \cite{Pl-Pukh-Analele}.

\vskip3mm
\section{The relations between monogenic functions and partial differential equations.}
\vskip3mm

Consider the following linear partial differential equation with
constant coefficients:
\begin{equation}\label{dopolnenije----1}
\mathcal{L}_NU(x,y,z):=\sum\limits_{\alpha+\beta+\gamma=N}C_{\alpha,\beta,\gamma}\,
\frac{\partial^N U}
{\partial x^\alpha\,\partial y^\beta\,\partial z^\gamma}=0, \quad
C_{\alpha,\beta,\gamma}\in\mathbb{R}.
\end{equation}

If a function $\Phi(\zeta)$  is $N$-times 
differentiable in the sense of Gateaux in every point of
$\Omega_{\zeta}$, then 
$$\frac{\partial^{\alpha+\beta+\gamma}\Phi}
{\partial x^\alpha\,\partial y^\beta\,\partial z^\gamma}=
e_1^\alpha\, e_2^\beta\, e_3^\gamma\,\Phi^{(\alpha+\beta+\gamma)}(\zeta)=
e_2^\beta\, e_3^\gamma\,\Phi^{(N)}(\zeta).
$$
Therefore, due to the equality
\begin{equation}\label{dopolnenije----2-1}
\mathcal{L}_N\Phi(\zeta)=\Phi^{(N)}(\zeta)
\sum\limits_{\alpha+\beta+\gamma=N}C_{\alpha,\beta,\gamma}\,e_2^\beta
\,e_3^\gamma\,,
\end{equation}
every $N$-times differentiable in the sense of Gateaux in
$\Omega_{\zeta}$ function $\Phi$ 
satisfies the equation $\mathcal{L}_N\Phi(\zeta)=0$ everywhere in
$\Omega_{\zeta}$ if and only if
\begin{equation}\label{dopolnenije----2}
\sum\limits_{\alpha+\beta+\gamma=N}C_{\alpha,\beta,\gamma}\,e_2^\beta\,
e_3^\gamma=0\,.
\end{equation}
Accordingly, if the condition (\ref{dopolnenije----2}) is
satisfied, then the real-valued components ${\rm Re}\,U_k(x,y,z)$
and ${\rm Im}\,U_k(x,y,z)$ of the decomposition
(\ref{rozklad-Phi-v-bazysi})
are solutions of the equation (\ref{dopolnenije----1}).

In the case where $f_u(E_3)=\mathbb{C}$ for all $u=1,2, \ldots,
m$, it follows from Theorem \ref{teo_pro_naslidky2} that the
equality (\ref{dopolnenije----2-1}) holds for every monogenic
function $\Phi:\Omega_{\zeta}\rightarrow \mathbb{A}_n^m$.

Thus, 
to construct solutions of the equation (\ref{dopolnenije----1}) in
the form of components of monogenic
functions, 
we must to find 
a triad of linearly independent over the field $\mathbb{R}$
vectors $(\ref{e_1_e_2_e_3})$ satisfying the characteristic
equation (\ref{dopolnenije----2}) and to verify the condition:
$f_u(E_3)=\mathbb{C}$ for all $u=1,2,\ldots,m$. Then, the formula
(\ref{Teor--1}) gives a constructive description of all mentioned
monogenic functions.

In the next theorem, we assign a special class of equations
(\ref{dopolnenije----1}) for which $f_u(E_3)=\mathbb{C}$ for all
$u=1,2,\ldots,m$. Let us introduce the polynomial
 \begin{equation}\label{dopolnenije----51}
P(a,b):=\sum\limits_{\alpha+\beta+\gamma=N}C_{\alpha,\beta,\gamma}\,
a^\beta\,b^\gamma\,.
\end{equation}
\vskip3mm

\theor\label{teo_dopolnenije-dlja-uravn} {\it Suppose that 
there exist linearly independent over the field $\mathbb{R}$
vectors $e_1,e_2,e_3$ in $\mathbb{A}_n^m$ of the form
$(\ref{e_1_e_2_e_3})$ that satisfy the equality
$(\ref{dopolnenije----2})$. If $P(a,b)\neq0$ for all real $a$ and
$b$, then $f_u(E_3)=\mathbb{C}$ for all $u=1,2,\ldots,m$.}

\textbf{Proof.} Using the multiplication table of 
$\mathbb{A}_n^m$, we obtain the equalities
$$e_2^\beta=\sum\limits_{u=1}^ma_u^\beta\,I_u+\Psi_\mathcal{R}\,,\quad
e_3^\gamma=\sum\limits_{u=1}^mb_u^\gamma\,I_u+\Theta_\mathcal{R}\,,
$$
where $\Psi_\mathcal{R}\,,\Theta_\mathcal{R}\in\mathcal{R}$. Now
the equality
 (\ref{dopolnenije----2}) takes the form
\begin{equation}\label{dopolnenije----4}
\sum\limits_{\alpha+\beta+\gamma=N}C_{\alpha,\beta,\gamma}\Bigr(
\sum\limits_{u=1}^ma_u^\beta\,b_u^\gamma\,I_u+\widetilde{\Psi}_\mathcal{R}\Bigr)=0,
\end{equation}
where $\widetilde{\Psi}_\mathcal{R}\in\mathcal{R}$. Moreover, due
to the assumption that the vectors $e_1,e_2,e_3$ of the form
$(\ref{e_1_e_2_e_3})$ satisfy the equality
$(\ref{dopolnenije----2})$, there exist complex coefficients
$a_k,b_k$ for $k=1,2,\ldots,n$ that satisfy the equality
(\ref{dopolnenije----4}).

It follows from the equality (\ref{dopolnenije----4}) that
\begin{equation}\label{dopolnenije----5}
\sum\limits_{\alpha+\beta+\gamma=n}C_{\alpha,\beta,\gamma}\,a_u^\beta\,b_u^\gamma=0,
\qquad u=1,2,\ldots,m.
\end{equation}

Since $P(a,b)\neq0$ for all $a,b\in\mathbb{R}$, the equalities
(\ref{dopolnenije----5}) can be satisfied only if for each
  $u=1,2,\ldots,m$ at least one of the numbers $a_u$ or $b_u$ belongs to
 $\mathbb{C}\setminus\mathbb{R}$
 that implies the relation $f_u(E_3)=\mathbb{C}$ for all\, $u=1,2,\ldots,m$.
The theorem is proved.

Note that 
if $P(a,b)\neq0$ for all $a,b\in\mathbb{R}$, then 
 $C_{N,0,0}\neq0$ because otherwise
 $P(a,b)=0$ for $a=b=0$.

Since the function $P(a,b)$ is continuous on $\mathbb{R}^2$, the
condition $P(a,b)\neq0$ means either $P(a,b)>0$ or $P(a,b)<0$ for
all $a,b\in\mathbb{R}$. Therefore, it is obvious that for any
equation (\ref{dopolnenije----1}) of elliptic type, the condition
$P(a,b)\neq0$ is always satisfied for all $a,b\in\mathbb{R}$.
At the same time, there are equations 
 (\ref{dopolnenije----1}) for which $P(a,b)>0$ for all $a,b\in\mathbb{R}$, but
which are not elliptic. For example, such are the equations
$$\frac{\partial^3 u}{\partial x^3}+\frac{\partial^3 u}{\partial x\partial y^2}+
\frac{\partial^3 u}{\partial x\partial z^2}=0
\;\;\;\mbox{and}\;\;\; \frac{\partial^5 u}{\partial
x^5}+\frac{\partial^5 u}{\partial x^3\partial y^2}+
\frac{\partial^5 u}{\partial x\partial y^2\partial z^2}=0
$$
considered in $\mathbb{R}^3$, in particular.

\vskip5mm \textbf{Acknowledgements.} The author expresses a
gratitude to Professor S.~A.~Plaksa
 and Mr. R.~P.~Pukhtaievych for numerous discussions 
 and valuable advices.

\vskip10mm

\vskip7mm
Vitalii Shpakivskyi

 Department of Complex Analysis and Potential Theory

 Institute of Mathematics of the National Academy of Sciences of
Ukraine,

 3, Tereshchenkivs'ka st.

01601 Kyiv-4

 UKRAINE

 http://www.imath.kiev.ua/\~{}complex/

 \ e-mail: shpakivskyi@mail.ru,\, shpakivskyi@imath.kiev.ua

\end{document}